\documentclass[a4paper]{article}
\usepackage{hyperref}
\usepackage{amsmath}
\usepackage{amsfonts}
\usepackage{amsthm}
\usepackage{indentfirst}
\usepackage{geometry}
\usepackage{color}

\renewcommand{\theequation}{\thesection.\arabic{equation}}
\usepackage{indentfirst}
\date{}
\makeatletter
\@addtoreset{equation}{section}

\newtheorem{theorem}{Theorem}[section]

\begin{document}
	
	\begin{center}
		{\Large  Asymptotic behavior of an integral equation of Landau-Lifschitz type}
	\end{center}
	
	\vskip 1cm
	\begin{center}
		{\sc Tiantian Zhou}
	\end{center}
	
	\begin{center}
		School of Mathematical Sciences\\
		Nanjing Normal University\\
		Nanjing, 210023, China\\
		email:zhoutiantiannj@163.com
	\end{center}
	
	\vskip 5mm
	{\leftskip5mm\rightskip5mm \normalsize \noindent{\bf{Abstract}}
		In this paper, we study an integral equation of Landau-Lifschitz-type involving the Riesz potential
		$$
			u(x)=\ell+C_*\int_{\mathbb{R}^n}\frac{u(y)(|\nabla u(y)|^{2}+|u_k(y)|^2)-u_k(y)e_k}{|x-y|^{n-\alpha}}dy,
		$$
		where $\ell \in \mathbb{R}^k$ is a constant vector, $C_* \neq 0$ is a constant,
		$u\colon\mathbb{R}^n \to \mathbb{S}^{k-1}:=\{u=(u_1,u_2,\cdots,u_k) \in
		{{\mathbb{R}}^{k}};|u|=1\}$, $k \geq 1$, $n \geq 2$, and $0<\alpha<n$.
		We provide a decay rate at which $|u_k|$ converges to zero under some mild conditions.
		This argument still works for the Ginzburg-Landau-type integral equation.
		\par
		\noindent{\bf{Keywords}}: Asymptotic behaviour, integral equation,
		Riesz potential, Landau-Lifschitz equation, Ginzburg-Landau equation
		
		\par
		{\bf{MSC}} 45M05, 35B40, 35Q56, 35Q60, 82D40}

	\newtheorem{proposition}[theorem]{Proposition}
	
	\renewcommand{\theequation}{\thesection.\arabic{equation}}
	\catcode`@=11
	\@addtoreset{equation}{section}
	\catcode`@=12
	
	\section{Introduction}

	In 2001, Hang and Lin studied the equation
	\begin{equation}\label{tll}
		-\Delta u=u(|\nabla u|^2+|u_3|^2)-u_3e_3 \quad {\rm on} \ \ \mathbb{R}^2,
	\end{equation}
	where $u\colon\mathbb{R}^2 \to \mathbb{S}^2:=\{u=(u_1,u_2,u_3) \in
	{{\mathbb{R}}^{3}};|u|=1\}$ and $e_3=(0,0,1)$. This is the Euler--Lagrange
	equation of
	$$
	E^{LL} (u)=\frac{1}{2}\|\nabla u\|_{L^2(\mathbb{R}^2)}^2
	+\frac{1}{2}\|u_3\|_{L^2(\mathbb{R}^2)}^2.
	$$
	Hang and Lin obtained the asymptotic rates of $u$ (cf. Theorem 1.4 in \cite{HL}):
	{\it If $0 \leq u_3 \in L^2(\mathbb R^2)$ and $|\nabla u(x)| \to 0$ $(|x| \to \infty)$,
		then} $|u_3(x)| \leq Ce^{-|x|/16}$.
	
	The energy functional $E^{LL} (u)$ comes from the theory of liquid crystals
	(cf.~\cite{AS}), the statics of planar ferromagnets and antiferromagnets
	(cf.~\cite{KP1,KP2,PS}) and the Schr\"odinger mapping equations
	(cf.~\cite{LW,WY}). In addition, according to Theorem~1.1 in \cite{GH},
	\eqref{tll} is related to the static Landau--Lifschitz equation which describes
	the dispersion of the magnetic permeability in ferromagnetic bodies
	(cf.~\cite{deL}). Thus, we call \eqref{tll} the static
	Landau--Lifschitz-type equation. In addition, We call a solution $u$ of \eqref{tll}
	the finite energy solution if $\nabla u \in L^2(\mathbb R^2)$, and call a solution $u$ of \eqref{tll}
	the finite penalization solution if $u_3 \in L^2(\mathbb R^2)$.
	Note that Remark 6.1 in \cite{HL} shows a Liouville theorem
	if $u$ is a finite energy solution. This result is stronger than the exponential decay
	property of the finite penalization solution.

	In this paper, we investigate the integral equation involving the Riesz potential
	\begin{equation}\label{1.1}
		u(x)=\ell+C_*\int_{\mathbb{R}^n}\frac{H(u(y),\nabla u(y))}{|x-y|^{n-\alpha}}dy,
	\end{equation}
	where $\ell \in \mathbb{R}^k$ is a constant vector, $C_* \neq 0$ is a constant, and
	$$
	H(u,\nabla u)=u(|\nabla u|^{2}+|u_k|^2)-u_ke_k, \ \ e_k=(0,\cdots,0,1).
	$$
	Here $u\colon\mathbb{R}^n \to \mathbb{S}^{k-1}:=\{u=(u_1,u_2,\cdots,u_k) \in
	{{\mathbb{R}}^{k}};|u|=1\}$, $k \geq 1$, $n \geq 2$, and $0<\alpha<n$. In view
	of \eqref{tll}, here \eqref{1.1} is also called the integral equation of the
	Landau--Lifschitz-type. If $u$ is a finite energy solution of \eqref{1.1},
	Liouville theorem holds true (cf. \cite{CL1,rjmp}). Here
	we are concerned with the asymptotic behavior of $|u_k|$ at infinity.
	
	\begin{theorem} \label{th1.1}
		Let $u \in C^{\infty}(\mathbb{R}^n,\mathbb{S}^{k-1})$ solve (\ref{1.1})
		with $n \geq 2$, $k \geq 1$ and $\alpha \in (0,n)$. Assume
		\begin{equation}\label{ti}
			u_k \in L^s(\mathbb{R}^n) \quad {\rm for \ some} \ s \in (1,n/\alpha),
		\end{equation}
		\begin{equation}\label{tii}
			\nabla u \in L^\infty(\mathbb{R}^n).
		\end{equation}
		Then $\ell_k=0$. Furthermore, if $u_k$ is radially symmetric and decreasing about
		the origin, then
		\begin{equation}\label{tiii}
			|u_k(x)| =O(|x|^{\alpha-n/s}) \quad (|x| \to \infty).
		\end{equation}
	\end{theorem}
	
	Different from Theorem 1.4 in \cite{HL}, we replace the finite penalization condition
	by \eqref{ti} here. In addition,
	to handle $H(u,\nabla u)$, it is necessary to impose restrictions on $|\nabla u|^2$.
	In \cite{HL}, the authors proved that the least energy solutions $u$ of
	\eqref{tll} satisfy \eqref{tii}. For other solutions $u$ of
	\eqref{tll}, De Laire \cite{deL} proved \eqref{tii} if $u_k \in L^2(\mathbb{R}^n)$.
	Therefore, if $s > 2$ in \eqref{ti}, the assumption \eqref{tii} is natural.

	The Landau-Lifschitz equation can be used to study the partial regularity of harmonic maps (cf. \cite{AS,HL}).
	The results in \cite{BBH,BMR} show that the Ginzburg-Landau equation has the same effect.
	
	In 1994, Shafrir \cite{Shafrir} studied the following Ginzburg-Landau equation
	\begin{equation}\label{1.5}
		-\Delta u=u(1-|u|^2) \quad {\rm on} \ \ \mathbb{R}^2.
	\end{equation}
	Here $u\colon\mathbb{R}^2 \to \mathbb{R}^2$ is a vector-valued function. This
	is the Euler--Lagrange equation for the Ginzburg--Landau energy
	$$
	E_{GL}(u)=\frac{1}{2}\|\nabla
	u\|_{L^2(\mathbb{R}^2)}^2+\frac{1}{4}\|1-|u|^2\|_{L^2(\mathbb{R}^2)}^2.
	$$
	Shafrir obtained the asymptotic rate of the finite penalization solution $u$
	as $|x| \to \infty$ (cf.~Theorem 1 in \cite{Shafrir}):
	$$
	1-|u|^2=d^2|x|^{-2}+o(|x|^{-2}) \quad (|x| \to \infty).
	$$
	where $d$ is the topological degree of $u(x)$ for large $|x|$.
	Theorem 2 in \cite{BMR} shows that the Liouville theorem still holds
	if $u$ is a finite energy solution of \eqref{1.5}.
	
	In 2005, papers \cite{MR,TZ} are concerned with fractional order Ginzburg--Landau equation:
	\begin{equation}\label{qw3}
		(-\Delta)^{\alpha/2}u=u(1-|u|^2) \quad {\rm on} \ \mathbb R^n.
	\end{equation}
	It describes the phase transition phenomena in fractal media.
	In 2016, Ma (cf.~\cite{Ma}) obtained the boundedness of finite penalization solutions of \eqref{qw3}.
	Those equations are also associated with the study of the fractional
	order harmonic mappings (cf.~\cite{MS}). Clearly, \eqref{qw3} is related to
	\begin{equation}\label{qw6}
		u(x)=\ell+C_*\int_{\mathbb R^n}\frac{u(y)(1-|u(y)|^2)}{|x-y|^{n-\alpha}}dy,
	\end{equation}
	Papers \cite{LX,LCL} gives the Liouville theorem
	for the finite energy solution of \eqref{qw6}.
	
	In 2021, Chen and Lei \cite{CL} studied
	the Allen--Cahn-type integral equations involving the Riesz potentials
	$$
	u(x)=\ell+C\int_{\mathbb R^n}\frac{H(u(y))dy}{|x-y|^{n-\alpha}}.
	$$
	Here $H(u)=u(1-|u|^2)|1-|u|^2|^{p-2}$ with $p>2$.
	The asymptotic rate of uniformly continuous and integrable solutions was given:
	$$
	|1-|u(x)|^2|=O(|x|^{\alpha-n}) \quad (|x| \to \infty).
	$$
	Afterwards, they proved the Liouville theorem for the finite energy solution (cf. \cite{CL2}).
	
	In this paper, we also consider the asymptotic behavior of integrable solutions
	of \eqref{qw6}.
	
	\begin{theorem} \label{th1.2}
		Let $u \in C^{\infty}(\mathbb{R}^n,\mathbb{R}^{k})$ solve (\ref{qw6})
		with $n \geq 2$, $k \geq 1$ and $\alpha \in (0,n)$. Assume
		\begin{equation}\label{tyi}
			1-|u|^2 \in L^s(\mathbb{R}^n) \quad {\rm for \ some} \ s \in (1,n/\alpha),
		\end{equation}
		and $1-|u|^2$ is radially symmetric and decreasing about
		the origin, then $u\to \ell$ as $|x|\to \infty$ with $|\ell|=1$, and
		\begin{equation}\label{tyiii}
			|1-|u(x)|^2| =O(|x|^{\alpha-n/s}) \quad (|x| \to \infty).
		\end{equation}
	\end{theorem}
	
	\section{Proof of Theorem \ref{th1.1}.}
	
	In this section, we use the Wolff potential and the
	Hardy--Littlewood--Sobolev inequality to investigate the value of $\ell_k$.
	
	{\it Step 1. We claim that $u_k(x) \to 0$ when $|x| \to \infty$ as long as \eqref{ti} and \eqref{tii} hold.}
	
	Otherwise, we can find $\epsilon_0>0$ and $|x_j| \to \infty$ such that $
	|u_k(x_j)| \geq 2\epsilon_0. $ Since $\nabla u \in L^\infty(\mathbb{R}^n)$,  it
	follows that there exists a $\eta>0$ such that $|u_k(x)-u_k(y)|<\epsilon_0$ for
	$|x-y|<\eta$. Therefore, $|u_k(x)|>\epsilon_0$ for $x \in B_\eta(x_j)$. Thus,
	\begin{equation}\label{2.1}
		\int_{B_\eta(x_j)}|u_k(x)|^s dx>\epsilon_0^s |B_\eta|.
	\end{equation}
	On the other hand, since $u_k \in L^s(\mathbb{R}^n)$, we have
	$$
	\int_{B_\eta(x_j)}|u_k(x)|^s dx
	\leq\int_{\{|x|>|x_j|-\eta\}}|u_k(x)|^s dx\to 0,
	$$ 
	This contradicts \eqref{2.1}.

	{\it Step 2. We claim} $\ell_k=0$.
	
	Write
	\begin{equation}\label{v-def}
		v(x)=\frac{1}{n-\alpha}
		\int_{\mathbb{R}^{n}}\frac{(|\nabla u(y)|^2+|u'(y)|^2)|u_k(y)|}{|x-y|^{n-\alpha}}dy.
	\end{equation}
	By exchanging the integral variants, we see that $v$ has the form involving the Wolff potential.
	Namely,
	$$
	v(x)=\int_0^{\infty}
	\frac{\int_{B_t(x)}[(|\nabla u(y)|^2+|u'(y)|^2)|u_k(y)|]dy}{t^{n-\alpha}}
	\frac{dt}{t}.
	$$
	
	Take $x_0 \in \mathbb{R}^n$. By $|u_k| \leq 1$, $|u'| \leq 1$ with $u'=(u_1,u_2,\dots, u_{k-1})$ and
	\eqref{tii}, we know that, for any $\varepsilon>0$, there exists a
	$\delta \in (0,1)$ such that
	$$
	\int_0^{\delta}
	\frac{\int_{B_t(x_0)}[(|\nabla u(y)|^2+|u'(y)|^2)|u_k(y)|]dy}{t^{n-\alpha}}
	\frac{dt}{t} \leq C
	\int_0^{\delta} t^{\alpha} \frac{dt}{t}
	<\varepsilon.
	$$
	If $|x-x_0|<\delta$, then $B_t(x_0) \subset B_{t+\delta}(x)$. Therefore,
	$$\begin{array}{ll}
		&\displaystyle\int_{\delta}^{\infty}
		\frac{\int_{B_t(x_0)} [(|\nabla u(y)|^2+|u'(y)|^2)|u_k(y)|]dy}{t^{n-\alpha}}
		\frac{dt}{t}\\
		\leq & C\displaystyle\int_{\delta}^{\infty}
		\frac{\int_{B_{t+\delta}(x)}[(|\nabla u(y)|^2+|u'(y)|^2)|u_k(y)|]dy}{(t+\delta)^{n-\alpha}}
		\left(\frac{t+\delta}{t}\right)^{n-\alpha+1} \frac{d(t+\delta)}{t+\delta}\\
		\leq & C2^{n-\alpha+1}\displaystyle\int_0^{\infty}
		\frac{\int_{B_t(x)}[(|\nabla u(y)|^2+|u'(y)|^2)|u_k(y)|]dy}{t^{n-\alpha}}
		\frac{dt}{t}\leq Cv(x).
	\end{array}
	$$
	Combining two estimates above, we obtain
	$$
	v(x_0)<\varepsilon+Cv(x) \quad {\rm when} \quad |x-x_0| <\delta.
	$$
	Thus, for any $t>1$,
	\begin{equation}\label{4.151}
		v^{t}(x_0)=|B_\delta(x_0)|^{-1}\displaystyle\int_{B_\delta(x_0)}v^{t}(x_0)dx
		\leq C\varepsilon^{t}
		+C|B_\delta(x_0)|^{-1}\displaystyle\int_{B_\delta(x_0)}v^{t}(x)dx.
	\end{equation}
	
	Applying the Hardy--Littlewood--Sobolev inequality to (\ref{v-def})
	and using $|u'| \leq 1$ and $\nabla u \in L^\infty(\mathbb{R}^n)$, we obtain
	\begin{equation}\label{qw1}
		\|v\|_{L^t(\mathbb{R}^n)} \leq
		C\|u_k\|_{L^{\frac{nt}{n+t\alpha}}(\mathbb{R}^n)}.
	\end{equation}
	Take $t=\frac{ns}{n-s\alpha}$. Then \eqref{qw1} implies $v \in L^{t}(\mathbb{R}^n)$.
	Therefore,
	$$
	\int_{B_\delta(x_0)}v^{t}(x)dx \to 0 \quad (|x_0| \to \infty)
	$$
	Inserting this result into (\ref{4.151}) with $t=\frac{ns}{n-s\alpha}$, we obtain
	\begin{equation}\label{rest}
		\lim_{|x| \to \infty}v(x)=0.
	\end{equation}

	It follows from (\ref{1.1}) that
	$$
	u_k(x)=\ell_k+C_*\int_{\mathbb{R}^n}\frac{[|\nabla u(y)|^2-|u'(y)|^2]u_k(y)}{|x-y|^{n-\alpha}}dy.
	$$
	Let $|x| \to \infty$. By \eqref{rest} and the result in Step 1, we see $\ell_k=0$.
	
	{\it Step 3. If $u_k$ is radially symmetric and decreasing about the origin.
		we claim that \eqref{tiii} holds true.}
	
	Denote $|u_k|$ by $V(x)$. In view of (\ref{tii}) and $|u'| \leq 1$, it follows
	\begin{equation*}
		V(x)\leq C\int_{{\mathbb{R}}^n}\frac{V(y)}{|x-y|^{n-\alpha}}dy.
	\end{equation*}
	Write
	$$
	T(\Omega):=\int_{\Omega}\frac{V(y)}{|x-y|^{n-\alpha}}dy.
	$$

	When $|y| \leq 3|x|/2$ and $|x-y| \geq \frac{|x|}{2}$,
	by the H\"older inequality and \eqref{ti}, we obtain
	\begin{equation}\label{t21}
		T\left(B(0,\frac{3|x|}{2}) \setminus B(x,\frac{|x|}{2})\right)
		\leq \frac{C}{|x|^{n-\alpha}}\|V\|_{L^s({\mathbb{R}}^n)}
		\left|B(0,\frac{3|x|}{2})\right|^{1-\frac{1}{s}}
		=\frac{C}{|x|^{\frac{n}{s}-\alpha}}.
	\end{equation}

	When $|y| \geq 3|x|/2$, there holds
	$
	|x-y| \geq |y|-|x| \geq |y|/3.
	$
	Therefore, by the H\"older inequality and \eqref{ti}, we obtain
	\begin{equation}\label{t22}
		T\left({\mathbb{R}}^n \setminus B(0,\frac{3|x|}{2})\right)
		\leq C\|V\|_{L^s({\mathbb{R}}^n)}
		\left(\int_{{\mathbb{R}}^n \setminus B(0,\frac{3|x|}{2})}
		|y|^{\frac{(\alpha-n)s}{s-1}}dy\right)^{1-\frac{1}{s}}
		\leq \frac{C}{|x|^{\frac{n}{s}-\alpha}}.
	\end{equation}

	Since $V$ is radially symmetric about the origin and decreasing, it holds that
	\begin{equation}\label{qw2}
		T\left(B(x,\frac{|x|}{2})\right) \leq CV(\frac{x}{2})|x|^\alpha.
	\end{equation}
	In addition, \eqref{ti} implies
	$$
	V^s(\frac{x}{2}) \left|B(0,\frac{|x|}{2})\right| \leq
	\int_{B(0,\frac{|x|}{2})} V^s(y)dy <\infty.
	$$
	This leads to $V(x/2) \leq C|x|^{-n/s}$. Inserting it into \eqref{qw2},
	we obtain
	$$
	T\left(B(x,\frac{|x|}{2})\right) \leq
	\frac{C}{|x|^{\frac{n}{s}-\alpha}}.
	$$
	Combining this result with (\ref{t21}) and \eqref{t22}, we obtain
	$$
	\int_{{\mathbb{R}}^n}\frac{V(y)}{|x-y|^{n-\alpha}}dy\leq 
	C\|V\|_{L^s(\mathbb{R}^n)}|x|^{\alpha-n/s}.
	$$ 
	This means that \eqref{tiii} holds.
	
	\section{Proof of Theorem \ref{th1.2}.}
	
	Denote $|1-|u|^2|$ by $V$. Since $V\in C^\infty(\mathbb{R}^n)$ is decreasing, $V\in L^\infty(\mathbb{R}^n)$, and hence  
	\begin{equation}\label{as3}
		u\in L^\infty(\mathbb{R}^n).
	\end{equation}

	Write
	\begin{equation*}
		v(x)=\frac{1}{n-\alpha}
		\int_{\mathbb{R}^{n}}\frac{|1-|u(y)|^2||u(y)|}{|x-y|^{n-\alpha}}dy.
	\end{equation*}
	By the same derivation of \eqref{rest}, using \eqref{as3}
	we also have $v(x) \to 0$ when $|x| \to \infty$. Combining this with \eqref{qw6} yields
	\begin{equation}\label{as1}
		u(x) \to \ell \quad (|x| \to \infty).
	\end{equation}

	We claim
	\begin{equation}\label{as2}
		|\ell|=1.
	\end{equation}
	Otherwise, from \eqref{as1} it follows
	$1-|u|^2 \to 1-|\ell|^2$ when $|x| \to \infty$. Thus, we can find $R>0$ such that
	$V(x) \geq |1-|\ell|^2|/2$ as $|x| >R$. Therefore,
	$$
	\int_{\mathbb R^n \setminus B_R(0)} V^s(x)dx
	\geq \left(\frac{|1-|\ell|^2|}{2}\right)^s \int_{\mathbb R^n \setminus B_R(0)} dx
	=\infty.
	$$
	It contradicts with \eqref{tyi}.
	
	By \eqref{as3} and \eqref{as2}, from \eqref{qw6} we deduce that
	$$
	V(x) =|(u+\ell)(u-\ell)| \leq |u-\ell||u+\ell|\leq C\int_{\mathbb R^n} \frac{V(y)dy}{|x-y|^{n-\alpha}}.
	$$
	By the same argument of Step 3 in the proof of Theorem \ref{th1.1}, we also obtain \eqref{tyiii}.
	
	\paragraph{Data availability.} Not applicable to this article as no datasets
	were generated or analysed during the current study.
	
	\paragraph{Disclosure statement.} This work does not have any conflicts of interest.

\end{document}